\documentclass[11pt]{amsart}
\usepackage{amsmath}
\usepackage{amsxtra}
\usepackage{amscd}
\usepackage{amsthm}
\usepackage{amsfonts}
\usepackage{amssymb}
\usepackage{eucal}
\newcommand{\bra}[1]{\langle #1 |}        
\newcommand{\ket}[1]{{| #1 \rangle}}      

\newcommand{\bea}{\begin{eqnarray}}
\newcommand{\ena}{\end{eqnarray}}
\newcommand{\be}{\begin{eqnarray*}}
\newcommand{\en}{\end{eqnarray*}}
\def\bel{\begin{eqnarray}}
\def\enl{\end{eqnarray}}

\newcommand{\C}{{\mathbb C}}

\newcommand{\Z}{{\mathbb Z}}

\newcommand{\al}{{\alpha}}

\newcommand{\la}{{\lambda}}
\newcommand{\mc}{\mathcal}

\numberwithin{equation}{section}

\numberwithin{equation}{section}
\newtheorem{thm}{Theorem}[section]
\newtheorem{prop}[thm]{Proposition}
\newtheorem{lem}[thm]{Lemma}
\newtheorem{cor}[thm]{Corollary}
\theoremstyle{remark}
\newtheorem{rem}[thm]{Remark}

\newcommand{\nc}{\newcommand}

\newcommand{\Id}{{\mathop {\rm id}}}

\newcommand{\hk}{\hookrightarrow}

\newcommand{\T}{\otimes}

\renewcommand{\deg}{\mathop{\rm deg}}

\renewcommand{\triangle}{\Delta}

\newcommand{\W}{\mathcal{W}}

\newcommand{\F}{\mathcal{F}}

\newcommand{\gli}{\mathfrak{gl}_\infty}
\newcommand{\hE}{\mathcal E}
\nc{\bc}{{\bf c}}
\nc{\bu}{{\bf u}}
\nc{\ba}{{\bf a}}

\begin{document}
\title[Semi-infinite construction of representations]
{Quantum continuous $\gli$: Semi-infinite construction of representations}

\author{B. Feigin, E. Feigin, M. Jimbo, T. Miwa and E. Mukhin}
\address{BF: Landau Institute for Theoretical Physics,
Russia, Chernogolovka, 142432, prosp. Akademika Semenova, 1a,   \newline
Higher School of Economics, Russia, Moscow, 101000,  Myasnitskaya ul., 20 and
\newline
Independent University of Moscow, Russia, Moscow, 119002,
Bol'shoi Vlas'evski per., 11}
\email{bfeigin@gmail.com}
\address{EF:
Tamm Department of Theoretical Physics, Lebedev Physics Institute, Russia, Moscow, 119991,
Leninski pr., 53 and\newline
French-Russian Poncelet Laboratory, Independent University of Moscow, Moscow, Russia }
\email{evgfeig@gmail.com}
\address{MJ: Department of Mathematics,
Rikkyo University, Toshima-ku, Tokyo 171-8501, Japan}
\email{jimbomm@rikkyo.ac.jp}
\address{TM: Department of Mathematics,
Graduate School of Science,
Kyoto University, Kyoto 606-8502,
Japan}\email{tmiwa@kje.biglobe.ne.jp}
\address{EM: Department of Mathematics,
Indiana University-Purdue University-Indianapolis,
402 N.Blackford St., LD 270,
Indianapolis, IN 46202}\email{mukhin@math.iupui.edu}

\begin{abstract}
We begin a study of the representation theory of quantum continuous
$\gli$, which we denote by $\hE$. This algebra depends on two parameters and
is a deformed version of the
enveloping algebra of the Lie algebra of difference operators acting on the
space of Laurent polynomials in one variable.
Fundamental representations of $\hE$
are labeled by a continuous parameter $u\in\C$. The representation
theory of $\hE$ has many properties
familiar from the representation theory of $\gli$: vector representations,
Fock modules, semi-infinite constructions of modules.
Using tensor products of vector representations,
we construct surjective homomorphisms from $\hE$ to spherical double affine Hecke algebras $S\ddot H_N$ for all $N$.
A key step in this construction is an identification
of a natural bases of the tensor products of vector representations
with Macdonald polynomials.
We also show that one of the Fock representations is isomorphic to the module constructed earlier
by means of the $K$-theory of Hilbert schemes.
\end{abstract}

\maketitle

\section{Introduction}
In this paper we begin to study the representation theory of an algebra $\hE$, which we call
the quantum continuous $\gli$.
This algebra is a
deformation of the universal enveloping algebra of the Lie algebra
of the $q$-difference operators in one variable.
Its representation theory
has a lot in common with that of the usual $\gli$ with a central extension:
vector representations,
fundamental representations, semi-infinite constructions of modules. Still there
is an important new feature of $\hE$: fundamental representations of this algebra
are labeled by a continuous parameter $u$. This makes the representation theory
of $\hE$ very rich and interesting. We give some details below.

The algebra $\hE$ is defined in terms of generators and relations. The generators
are denoted by $e_i$, $f_i$ ($i\in\Z$) and $\psi^+_j$, $\psi^-_{-j}$ ($j\ge 0$).
The elements $\psi_0^\pm$ are central and invertible. The relations between
generators depend symmetrically  on three parameters $q_1,q_2,q_3$ which are assumed to satisfy
$q_1q_2q_3=1$. These relations are given explicitly in terms of generating series
(see Section \ref{definition}). For example, let $e(z)=\sum_{i\in\Z} e_i z^{-i}$.
Then the following relation holds in $\hE$:
\begin{equation}\label{ee}
g(z,w)e(z)e(w)=-g(w,z)e(w)e(z), \quad g(z,w)=(z-q_1w)(z-q_2w)(z-q_3w).
\end{equation}
This relation appears in different contexts (see \cite{FO}, \cite{Kap}, \cite{FT},
\cite{SV2}).
In terms of components \eqref{ee} is equivalent to the set of relations
labeled by integers $m,n\in\Z$:
\begin{multline*}
e_{n+3}e_m-(q_1+q_2+q_3)e_{n+2}e_{m+1}
+(q_1q_2+q_1q_3+q_2q_3)e_{n+1}e_{m+2}-e_ne_{m+3}
\\
=
-e_{m+3}e_n + (q_1+q_2+q_3)e_{m+2}e_{n+1}
- (q_1q_2+q_1q_3+q_2q_3)e_{m+1}e_{n+2} + e_me_{n+3}.
\end{multline*}

As we have mentioned above, relation \eqref{ee} has different origins.
Let us explain the one important for us.
Fix a parameter $q$
and consider the associative algebra
$A=\C[Z,Z^{-1},D,D^{-1}]$ with $DZ=qZD$.
The algebra $A$ acts on the space $\C[z,z^{-1}]$ by
$Z(f)(z)=z f(z)$, $(Df)(z)=f(qz)$.
Thus $A$ can be identified with the algebra of $q$-difference
operators, and can be thought of as an algebra of special
infinite matrices.
The algebra $A$ admits representations with a continuous parameter $u$ on the space of delta
functions $\bigoplus_{i\in\Z} \delta(q^iu/z)$ through the same action on a vector $f(z)$ in this space.
Thus $A$ may be called a continuous $\gli$.

Consider the elements $\bar e_i=Z^iD\in A$.
It is easy to check that these elements satisfy the relations
\begin{multline*}
\bar e_{n+3}\bar e_m-(1+q+q^{-1})\bar e_{n+2}\bar e_{m+1}
+(1+q+q^{-1})\bar e_{n+1}\bar e_{m+2}-\bar e_n\bar e_{m+3}
\\
=
-\bar e_{m+3}\bar e_n + (1+q+q^{-1})\bar e_{m+2}\bar e_{n+1}
- (1+q+q^{-1})\bar e_{m+1}\bar e_{n+2} + \bar e_m\bar e_{n+3}.
\end{multline*}
Thus the relation \eqref{ee} is a quantization of the relations above.
In fact all other relations of $\hE$
(Section \ref{definition}) can be obtained in a similar way.
Moreover, there exists a Poisson structure on $A$ (considered as a Lie algebra)
such that the usual quantization technique, applied
to the universal enveloping algebra of $A$, gives $\hE$. We do not discuss this construction
in this paper and
will return to it elsewhere.

We recall that in \cite{DI} the authors constructed a class of quantum algebras generalizing
quantum affine algebras. A particular example of their construction is
the algebra $\hE'$, which differs from $\hE$ only by the absence of
the cubic relations (see \eqref{rel6} below)
$$
[e_0,[e_1,e_{-1}]]=0, \quad [f_0,[f_1,f_{-1}]]=0.
$$
We call these relations Serre relations for $\hE$.
The algebra $\hE'$ was also considered in \cite{FHHSY} and \cite{FT} and was
called there the Ding-Iohara algebra.
We note however that the Serre relations are important from the point of view of the
representation theory of
$\hE$ and of the structure theory as well. We explain the reasons below.

We recall that in \cite{SV2}, \cite{FT} the equivariant localized $K$-theory of
Hilbert schemes $H_n$ of $n$ points of $\C^2$ was studied.
In particular it was shown that the direct sum
$\F=\bigoplus_{n\ge 0} K(H_n)$ is isomorphic to
the space of  symmetric polynomials in infinite number of variables
and carries the structure of $\hE'$ module.
We note however that this action factors through the surjection
$\hE'\to\hE$ and therefore $\F$ has a natural structure of $\hE$-module.
The space  $\F$
has a natural basis labeled by
fixed points of the action of the torus on $H_n$ and
elements of this basis can be identified with the
Macdonald polynomials. In addition the action of the generators $e_i$ and $f_i$
is given by Pieri-like formulas.
In this paper we observe that
the representation $\F$ can be constructed
by means of a version of the semi-infinite
wedge construction.

Recall that the main building block of the semi-infinite construction for $\gli$
is its vector representation.
We start with considering $\hE$-modules
$V(u)$ ($u\in\C$) which is spanned by the vectors $[u]_i$ ($i\in\Z$). They
play the role of the vector representation.
The usual $\gli$ has only one vector
representation, but $\hE$ naturally has
a continuous family of such representations.
The algebra $\hE$ is endowed with a structure of ``comultiplication" (see \cite{DI}).
Strictly speaking, this ``comultiplication'' does not define
a structure of $\hE$-module on an arbitrary tensor product $V\otimes W$ of $\hE$-modules,
because some convergence conditions need to be satisfied (see Section \ref{definition}
for details).
We show that the tensor product
$V(u_1)\otimes\cdots\otimes V(u_N)$
is well-defined for general values of  $u_1,\dots,u_N$.
We are mainly interested in the case when the parameters $u_i$ form a geometric
progression.
We show that the tensor product
$$V(u)\otimes V(u q_2^{-1})\otimes \cdots\otimes V(u q_2^{-N+1})$$
has a subrepresentation
$W^N(u)$ spanned by the set of vectors
$[u]_{i_1}\otimes [u q_2^{-1}]_{i_2}\otimes
\cdots\otimes [u q_2^{-N+1}]_{i_N}$ with
$i_1> i_2>\cdots> i_N$.
The $\hE$-modules $W^N(u)$
are analogues of the exterior powers of the vector
representation for $\gli$.
We construct the structure of $\hE$-module on the limit
$N\to\infty$ of $W^N(u)$,
thus obtaining an analogue of the space of
semi-infinite forms. We denote this representation by $\F(u)$
and call it the Fock representation.

The space
$W^N(u)$
can be identified with the space of symmetric polynomials
in $N$ variables.
We recall that the space
$\C[x_1^{\pm1},\cdots,x_N^{\pm1}]^{S_N}$
has a natural structure of faithful representation of the spherical double affine Hecke algebra
$S\ddot H_N$. We show that the image of $\hE$ coincides with spherical DAHA and thus
obtain a surjective homomorphism $\hE\to S\ddot H_N$ for any $N$. We recall
that in \cite{SV1}, \cite{SV2} the spherical DAHA of type $GL_\infty$
was constructed as a projective limit of $S\ddot H_N$. It is
natural to expect that our $\hE$ is isomorphic to $\lim_{N\to\infty} S\ddot H_N$
(we plan to discuss this elsewhere).

Because of the homomorphisms $\hE\to S\ddot H_N$ any
$S\ddot H_N$-module gives us a representation of $\hE$.
 Consider now the resonance case
$q_1^{1-r}q_3^{k+1}=1$, $k>0$, $r>1$ .
In this case, the representation
of
$S\ddot H_N$
on
$\C[x_1^{\pm1},\cdots,x_N^{\pm1}]^{S_N}$
has a subrepresentation
$W^{k,r,N}
\subset \C[x_1^{\pm1},\cdots,x_N^{\pm1}]^{S_N}$ defined by
$$
W^{k,r,N}=\{f(x_1,\cdots,x_N)\ |\ f({\boldsymbol x})=0 \text{ if } {\boldsymbol x}
\text{ satisfies the wheel condition } \},
$$
where where ${\boldsymbol x}=(x_1,\cdots,x_N)$ is
said to satisfy the wheel condition if
$$
x_i=x_1 q_3^{1-i} q_1^{s_1+\dots + s_{i-1}}, i=1,\dots,k+1,\ s_1,\dots,s_{k+1}\ge 0,\ s_1+\dots + s_{k+1}=r-1.
$$
In \cite{FJMM2} it is proved that
$W^{k,r,N}$
has a basis labeled by the so-called
$(k,r)$-admissible partitions, i.e. partitions $\la$ satisfying $\la_i-\la_{i+k}\ge r$ for all $i\ge 1$.
Each element of the basis is a Macdonald (Laurent) polynomial.
Thus we have an action of the algebra $\hE$
on the space of polynomials satisfying the wheel condition.
We construct a family of $\hE$-modules $W^{k,r,N}(u)$, $u\in\C$ such that
$W^{k,r,N}(1)\simeq W^{k,r,N}$.
We also construct the inductive limit $N\to\infty$
of the modules $W^{k,r,N}$ and endow it with a structure of the
$\hE$-module.
As a result, we construct a family of representations $W^{k,r}_{\bc}(u)$
of  $\hE$  whose bases are labeled by infinite $(k,r)$-admissible partitions
with certain stability property at infinity. The parameter $\bc=(c_1,\ldots,c_{k-1})$
($1\leq c_1\leq\cdots\leq c_{k-1}\leq r$) enters in the stability property.

Our paper is organized as follows. \\
In Section $2$ we give the definition of $\hE$.\\
In Section $3$ the vector representations and their tensor products are constructed.\\
In Section $4$ we work out the semi-infinite construction for general parameters $q_i$.\\
In Section $5$ we establish a link between the tensor products of representations of
$\hE$ and representations of $S\ddot H_N$.\\
In Section $6$ we consider the semi-infinite construction in
the resonance case $q_1^{1-r}q_3^{k+1}=1$.\\
In Section $7$ we discuss further properties of
the algebra $\hE$.

\section{Quantum continuous $\gli$}\label{definition}
In this section we introduce the algebra $\hE$ which we call the quantum continuous
$\gli$.

\subsection{Definition}
Let $q_1,q_2,q_3$ be complex numbers, satisfying $q_i\ne 1$ and
$q_1q_2q_3=1.$
Let
\be
g(z,w)=(z-q_1w)(z-q_2w)(z-q_3w).
\en
Let $\hE$ be an  associative algebra over $\C$ generated by the elements $e_i$, $f_i$
($i\in\Z$), $\psi^+_j$, $\psi^-_{-j}$ ($j> 0$) and $(\psi^\pm_0)^{\pm 1}$
with defining relations depending on parameters $q_1,q_2,q_3$.
(So strictly speaking, we have a family of algebras).
We use  generating series
$$
e(z)=\sum_{i\in\Z} e_iz^{-i},\quad  f(z)=\sum_{i\in\Z} f_iz^{-i},\quad
\psi^\pm(z)=\sum_{\pm i\ge 0}\psi^{\pm}_iz^{-i}.
$$
The defining relations in $\hE$ are
\begin{gather}
\label{rel1}
g(z,w)e(z)e(w)=-g(w,z)e(w)e(z),\qquad
g(w,z)f(z)f(w)=-g(z,w)f(w)f(z),\\
\label{rel2}
g(z,w)\psi^{\pm}(z)e(w)=-g(w,z)e(w)\psi^{\pm}(z),\quad
g(w,z)\psi^{\pm}(z)f(w)=-g(z,w)f(w)\psi^{\pm}(z),\\
\label{rel3}
[e(z), f(w)]=\frac{\delta(z/w)}{g(1,1)}(\psi^{+}(z)-\psi^{-}(z)),\\
\label{rel4}
[\psi^\pm_i,\psi^\pm_j]=0, \qquad [\psi^\pm_i,\psi^\mp_j]=0, \\
\label{rel5}
\psi^\pm_0(\psi^\pm_0)^{-1}=(\psi^\pm_0)^{-1}\psi^\pm_0=1,\\
\label{rel6}
[e_0,[e_1,e_{-1}]]=0,\quad [f_0,[f_1,f_{-1}]]=0.
\end{gather}
Here $\delta(z)=\sum_{n\in\Z} z^n$ is the delta-function.

\medskip
\begin{rem}
The form of relations \eqref{rel1}, \eqref{rel2}, \eqref{rel3}
is a convenient way of writing algebraic relations between generators.
Namely each relation is to be understood as generating functions for relations for Fourier coefficients
of the right and left hand sides. For example, the relation
\eqref{rel1} for $e(z)$ is equivalent to the following set of relations labeled by  pairs $(n,m)\in\Z^2$:
\begin{multline*}
e_{n+3}e_m - (q_1+q_2+q_3)e_{n+2}e_{m+1} + (q_1q_2+q_2q_3+q_3q_1)e_{n+1}e_{m+2} - e_ne_{m+3}\\=
-e_{m+3}e_n + (q_1+q_2+q_3)e_{m+2}e_{n+1} - (q_1q_2+q_2q_3+q_3q_1)e_{m+1}e_{n+2} + e_me_{n+3}
\end{multline*}
and  \eqref{rel3} simply means that
$$
g(1,1)[e_i,f_j]=
\begin{cases}
\psi^+_{i+j}, \text{ if } i+j>0,\\
-\psi^-_{i+j}, \text{ if } i+j<0,\\
\psi^+_0-\psi^-_0, \text{ if } i+j=0.
\end{cases}
$$
\end{rem}
\medskip

\begin{rem}\label{field}
The algebra $\hE$ can be considered as an algebra over one of the fields
of rational functions $\C(q_1,q_2)$, $\C(q_1,q_3)$, $\C(q_1,q_3)$.
This is equivalent to saying that
the parameters $q_i$ are ``general". However, $\hE$ is defined for
arbitrary (except $1$ and $0$) values of parameters. Also in Section \ref{secres}
we will consider the case when $q_i$ satisfy an algebraic relation
(different from $q_1q_2q_3=1$).
\end{rem}
\medskip

In what follows we call the algebra  $\hE$ the quantum continuous $\gli$.
\medskip
\begin{rem}\label{DI}
In \cite{DI} Ding and Iohara defined a class of algebras, which are analogues
of quantum affine algebras.
Apart from the cubic Serre relations \eqref{rel6},
the algebra $\hE$ is a particular case of their construction.
This algebra (without relations \eqref{rel6}) was also considered
in \cite{FT}, \cite{FHHSY}.
\end{rem}
\medskip

The following lemma is obvious.

\begin{lem}
The algebra $\hE$ is invariant under permutations of parameters
$q_1,q_2,q_3$.

The elements $\psi_0^\pm\in \hE$ are central.

There is an anti-involution of $\hE$ sending $e(z)$ to
$f(z)$, $f(z)$ to $e(z)$ and $\psi^{\pm}(z)$ to $\psi^{\mp}(z)$.

The algebra $\hE$ is graded by the lattice $\Z^2$. The degrees of generators are given by
\[
\deg e_i= (1,i),\qquad \deg f_i=(-1,i),\qquad \deg \psi^\pm_i= (0,i).
\]
\end{lem}

We say that an  $\hE$-module is of level $(l_+,l_-)$ if $\psi^\pm_0$ act on this representation by scalars $l_\pm$.

Let $\hE'$ be the algebra defined in the same way as $\hE$ without cubic relations
$\eqref{rel6}$ (see Remark \ref{DI}).
In \cite{DI} the formal (see the explanations below)
structure of the Hopf algebra on $\hE'$ was constructed.
In particular the comultiplication
is given by
\begin{gather}\label{tre}
\triangle e(z)= e(z)\T 1 + \psi^-(z)\T e(z),\\
\label{trf}
\triangle f(z)= f(z)\T \psi^+(z) + 1\T f(z),\\
\label{trpsi}
\triangle \psi^\pm(z)= \psi^\pm(z)\T \psi^\pm(z).
\end{gather}
We note that this ``definition" does not define a comultiplication in
the usual sense.
The right hand sides are not elements of $\hE\T\hE$ since they contain
infinite sums. Still for certain classes of modules the formulas
\eqref{tre}, \eqref{trf} and \eqref{trpsi} can be made precise.
So in what follows, when talking about the tensor products
$V_1\T\cdots\T V_N$ of $\hE$-modules, we construct the action of
the generators $e_i$, $f_i$ and $\psi^\pm_i$ explicitly
(based on the universal formulas \eqref{tre}, \eqref{trf}, \eqref{trpsi})
and check that they satisfy the relations of quantum continuous $\gli$.

We close this section with the following statement.
\begin{lem}\label{SERRE}
In $\hE'$ the cubic element $[e_0,[e_1,e_{-1}]]$ belongs to the kernel of ${\rm ad}\,f(z)$. Similarly,
$[f_0,[f_1,f_{-1}]]$ belongs to the kernel of ${\rm ad}\,e(z)$.
\end{lem}
The proof will be given elsewhere.
Using this lemma it is not difficult to prove the Serre relations in each representation we discuss
in this paper.

\section{The modules $V(u)$}
In this section we define vector representations $V(u)$ of $\hE$.
We also construct tensor products of vector representations and 
certain submodules inside tensor products.
\subsection{Vector representations}
For a parameter $u\in\C$ we consider
the space $V(u)$,  spanned by basis vectors $[u]_i$ ($i\in\Z$).
In the following lemma we define representations
of the quantum continuous $\gli$ depending on parameter $u$.
We call $V(u)$ a vector representation.
\begin{prop}\label{V_u}
The assignment
\begin{align*}
&(1-q_1)e(z)[u]_i=\delta(q_1^iu/z)[u]_{i+1}\,,
\\
&-(1-q_1^{-1})f(z)[u]_i=\delta(q_1^{i-1}u/z)[u]_{i-1}\,,
\\ &\psi^+(z)[u]_i=\frac{(1-q_1^iq_3u/z)(1-q_1^iq_2u/z)}{(1-q_1^iu/z)(1-q_1^{i-1}u/z)}[u]_i,\\
&\psi^-(z)[u]_i=\frac{(1-q_1^{-i}q_3^{-1}z/u)(1-q_1^{-i}q_2^{-1}z/u)}{(1-q_1^{-i}z/u)(1-q_1^{-i+1}z/u)}[u]_i,
\end{align*}
defines a structure of level $(1,1)$ $\hE$-module on $V(u)$.
\end{prop}

\medskip
\begin{rem}
An important feature of the representations $V(u)$ is that  $\psi^\pm(z)$ act
on $[u]_i$ via multiplication by the expansions at $z=\infty$ and $z=0$ of the function
\[
\frac{(1-q_1^{-i}q_3^{-1}z/u)(1-q_1^{-i}q_2^{-1}z/u)}{(1-q_1^{-i}z/u)(1-q_1^{-i+1}z/u)}.
\]
\end{rem}
\medskip

For the proof of Proposition \ref{V_u} we need a simple lemma.
We use the following notation: for a rational
function $\gamma(z)$ we denote by $\gamma^\pm(z)$ the expansions of $\gamma(z)$
at $z=\infty$ and $z=0$,
i.e. $\gamma^\pm(z)$ are Taylor series in $z^{\mp 1}$.
\begin{lem}\label{sp}
Let $\gamma(z)$ be a rational function regular at $z=0,\infty$ and with simple
poles. Then we have the formal series identity
\[
\gamma^+(z) - \gamma^-(z) = \sum_t \gamma^{(t)}\delta(z/z^{(t)}),
\]
where the sum runs over all poles $z^{(t)}$ of $\gamma(z)$ and
$\gamma^{(t)}=res_{z=z^{(t)}} \gamma(z)\frac{dz}{z}$.
\end{lem}

We now prove Proposition \ref{V_u}.
\begin{proof}
Since $e_m [u]_j=(1-q_1)^{-1}q_1^{jm}u^m[u]_{j+1}$, the relations \eqref{rel6} are obviously satisfied.
We show now that \eqref{rel1} and \eqref{rel3} hold, all other relations are proved
similarly. In what follows we often use the formula
\begin{equation}\label{delta}
\gamma(z)\delta(z/w)=\gamma(w)\delta(z/w).
\end{equation}

So let us show that
\[
g(z,w)e(z)e(w)=-g(w,z)e(w)e(z).
\]
In fact, we prove that both sides vanish on $V(u)$. By definition
\begin{align*}
(1-q_1)^2g(z,w)e(z)e(w)[u]_i
&=g(z,w)\delta(q_1^{i+1}u/z)\delta(q_1^i u/w)[u]_{i+2}\\
&=g(q_1^{i+1}u,q_1^iu)\delta(q_1^{i+1}u/z)\delta(q_1^i u/w)[u]_{i+2}\\
&=0.
\end{align*}
Similarly, $g(w,z)e(w)e(z)=0.$

Now we show that
\begin{equation}\label{ef}
[e(z), f(w)][u]_i=\frac{\delta(z/w)}{g(1,1)}(\psi^{+}(w)-\psi^{-}(z))[u]_i.
\end{equation}
The left hand side reads as
\begin{multline}\label{d-d}
\frac{q_1}{(1-q_1)^2}(\delta(q_1^{i-1}u/w)\delta(q_1^{i-1}u/z) -
\delta(q_1^i u/w)\delta(q_1^iu/z))[u]_i\\ =
\frac{q_1}{(1-q_1)^2}\delta(z/w)(\delta(q_1^{i-1}u/z) - \delta(q_1^i u/z)).
\end{multline}
The right hand side of \eqref{ef} equals to
\begin{equation}\label{frac}
\frac{\delta(z/w)}{g(1,1)}
\left(
\frac{1-q_1^iq_3u/z}{1-q_1^iu/z}\times\frac{1-q_1^iq_2u/z}{1-q_1^{i-1}u/z} -
\frac{1-q_1^{-i}q_3^{-1}z/u}{1-q_1^{-i}z/u}\times \frac{1-q_1^{-i}q_2^{-1}z/u}{1-q_1^{-i+1}z/u}
\right)[u]_i.
\end{equation}
Since the expression in the round brackets is of the form
$\gamma^+(z) - \gamma^-(z)$ for a rational function $\gamma(z)$, we can apply
Lemma \ref{sp}, which proves \eqref{ef}.
\end{proof}

We define the rational functions
\[
\gamma_{i,u}(z)=\frac{(1-q_3q_1^iu/z)(1-q_2q_1^iu/z)}{(1-q_1^{-1}q_1^iu/z)(1-q_1^iu/z)}.
\]
Then we have
$
\psi^\pm(z)[u]_i=\gamma_{i,u}^\pm(z)[u]_i.
$

\subsection{Tensor products}
Consider the tensor product of vector representations $V(u_1)\T\dots\T V(u_N)$.
We define the following generating series of operators on this space:
\begin{multline}\label{ereg}
(1-q_1)e(z) ([u_1]_{i_1}\T\dots\T [u_N]_{i_N})\\
=\sum_{s=1}^N
\left(\prod_{l=1}^{s-1} \gamma_{i_l,u_l}(q_1^{i_s}u_s)\right)
\delta(q_1^{i_s}u_s/z)
[u_1]_{i_1}\T\dots \T [u_{s-1}]_{i_{s-1}}\T  [u_s]_{i_s+1}\T [u_{s+1}]_{i_{s+1}}\T \dots\T  [u_N]_{i_N},
\end{multline}
\begin{multline}\label{freg}
-(1-q_1^{-1})f(z) ([u_1]_{i_1}\T\dots\T [u_N]_{i_N})\\
=\sum_{s=1}^N
\delta(q_1^{i_s-1}u_s/z)\left(\prod_{l=s+1}^N \gamma_{i_l,u_l}(q_1^{i_s-1}u_s)\right)
[u_1]_{i_1}\T\dots \T [u_{s-1}]_{i_{s-1}}\T  [u_s]_{i_s-1}\T [u_{s+1}]_{i_{s+1}}\T \dots\T  [u_N]_{i_N},
\end{multline}
\begin{equation}\label{psireg}
\psi^\pm(z) ([u_1]_{i_1}\T\dots\T [u_N]_{i_N})
=\psi^\pm(z)[u_1]_{i_1}\T\dots \T  \psi^\pm(z) [u_N]_{i_N}.
\end{equation}

The formulas above are read from the universal formulas \eqref{tre}, \eqref{trf},
\eqref{trpsi}. In fact, formula \eqref{tre} gives (formally) the action of $e(z)$
on the tensor product $V(u_1)\T\dots\T V(u_N)$:
\begin{equation}\label{formal}
e(z) = \sum_{s=1}^N \underbrace{\psi^-(z) \T \dots \psi^-(z)}_{s-1} \T e(z) \T \Id\T\dots\T\Id.
\end{equation}
By definition, $e(z)[u]_i=(1-q_1)^{-1}\delta(q_1^iu/z)[u]_{i+1}$ and $\psi^\pm(z)$ acts on
$[u]_i$ via multiplication by certain series. The product of delta function
with these series is in general not defined.
The series in question are expansions of given rational functions.
It is therefore natural to regularize \eqref{formal} by substituting  the support
$z=q_1^iu$ of the delta function into rational functions. We thus obtain
formula \eqref{ereg}. Similar arguments lead to \eqref{freg}.
Note however that the expression $\gamma_{i_l,u_l}(q_1^{i_s}u_s)$ is not defined
if the argument is a pole of $\gamma_{i_l,u_l}(z)$. Therefore, formulas
\eqref{ereg} and \eqref{freg} do not always produce well-defined operators.

For an element $\ba=(a_1,\dots,a_N)\in\Z^N$ let
$\bu_\ba\in \bigotimes_{s=1}^N V(u_s)$ denote the vector
$\bigotimes_{s=1}^N [u_s]_{a_s}$.
\begin{lem}\label{sub}
Let $A\subset \Z^N$ be a subset such that
\begin{itemize}
\item for all $\ba\in\Z^N$, $\ba'\in A$ the matrix coefficients
$\bra{\bu_{\ba'}} e(z) \ket{\bu_\ba}$ and
$\bra{\bu_{\ba'}} f(z) \ket{\bu_\ba}$ are well-defined,
\item for all $\ba\in A$, ${\bf b}\notin A$ the matrix coefficients
$\bra{\bu_{\bf b}} e(z) \ket{\bu_\ba}$ and $\bra{\bu_{\bf b}} f(z) \ket{\bu_\ba}$ vanish.
\end{itemize}
Then formulas \eqref{ereg}, \eqref{freg} and \eqref{psireg} define a structure
of $\hE$-module on $\mathrm{span}\{\bu_\ba\}_{\ba\in A}$.
\end{lem}
\begin{proof}
It suffices to show that the defining relations of $\hE$ are satisfied.
We check \eqref{rel3} and \eqref{rel6}. The rest can be checked similarly.

We start with relation \eqref{rel3}. Let us compute the matrix coefficient
\begin{equation}\label{aa}
\bra{\bu_{\ba'}} [e(z),f(w)]\ket{\bu_\ba}.
\end{equation}
We first show that it vanishes unless $\ba'=\ba$. Introduce the notation
$\ba \pm {\bf 1}_s=(...,a_s\pm 1, ...)$. Then clearly \eqref{aa} vanishes unless
$\ba'=\ba + {\bf 1}_s - {\bf 1}_t$ for some $s,t=1,\dots,N$.
If $s<t$, then \eqref{aa} vanishes because formulas \eqref{ereg} and \eqref{freg}
give identical expressions for
$\bra{\bu_{\ba'}} e(z)f(w)\ket{\bu_\ba}$ and for $\bra{\bu_{\ba'}} f(w)e(z)\ket{\bu_\ba}$.
Assume $s< t$. Then
from formulas \eqref{ereg} and \eqref{freg} we obtain that \eqref{aa} is equal to
\begin{multline*}
\left(\gamma_{a_t,u_t}(q_1^{a_s-1}u_s)\gamma_{a_s-1,u_s}(q_1^{a_t}u_t) -
\gamma_{a_s,u_s}(q_1^{a_t}u_t)\gamma_{a_t+1,u_t}(q_1^{a_s-1}u_s)\right)\\
\times\prod_{s\ne l<t} \gamma_{a_l,u_l}(q_1^{a_t}u_t)
\prod_{t\ne l>s} \gamma_{a_l,u_l}(q_1^{a_s-1}u_s)
\delta(q_1^{a_s-1}u_s/w)\delta(q_1^{a_t}u_t/z),
\end{multline*}
which vanishes thanks to a simple relation
\[
\gamma_{a_t,u_t}(q_1^{a_s-1}u_s)\gamma_{a_s-1,u_s}(q_1^{a_t}u_t) =
\gamma_{a_s,u_s}(q_1^{a_t}u_t)\gamma_{a_t+1,u_t}(q_1^{a_s-1}u_s).
\]
So we only have the terms with $s=t$ and thus $\ba\ne \ba'$ implies \eqref{aa} is zero.

Now assume $\ba=\ba'$. Then \eqref{aa} is equal to
\begin{multline}
\sum_{s=1}^N \delta(q_1^{a_s-1}u_s/w) \delta(q_1^{a_s-1}u_s/z)
\prod_{l\ne s} \gamma_{a_l,u_l}(q_1^{a_s-1}u_s) -
\sum_{s=1}^N \delta(q_1^{a_s}u_s/w) \delta(q_1^{a_s}u_s/z)
\prod_{l\ne s} \gamma_{a_l,u_l}(q_1^{a_s}u_s)\\
=\delta(z/w)
\left(
\sum_{s=1}^N \delta(q_1^{a_s-1}u_s/z)
\prod_{l\ne s} \gamma_{a_l,u_l}(q_1^{a_s-1}u_s) -
\sum_{s=1}^N \delta(q_1^{a_s}u_s/z)
\prod_{l\ne s} \gamma_{a_l,u_l}(q_1^{a_s}u_s)
\right).
\end{multline}

By definition of the action of $\psi^\pm(z)$ we have
$$
\langle \mathbf{u}_{\mathbf{a}}|\bigl(\psi^+(u)-\psi^-(u)\bigr)|\mathbf{u}_{\mathbf{a}}\rangle
$$
Assume for a moment that
\begin{equation}\label{qu}
q_1^i u_l \ne q_1^j u_m\quad \text{ unless }\quad i=j, l=m.
\end{equation}
Then
all poles of the function $\prod_{s=1}^N \gamma_{a_s,u_s}(z)$ are simple and
Lemma \ref{sp} proves
\begin{equation}\label{efpsi}
[e(z), f(w)]\bu_\ba = \frac{\delta(z/w)}{g(1,1)} (\psi^+(z) - \psi^-(z)) \bu_\ba.
\end{equation}
We note also that if relation \eqref{efpsi} holds for parameters satisfying
\eqref{qu}, then it holds for all values of parameters.

We now prove relation \eqref{rel6}. Let $E=[e_0,[e_1,e_{-1}]]$.
Let $\ba\pm{\bf 1}_j=(a_1,\ldots,a_j\pm 1,\ldots,a_N)$.
From Lemma \ref{SERRE} follows for all $1\leq i\leq j\leq N$ that
$$
\sum_{n=1}^N\bra{\bu_{\ba +{\bf 1}_i+{\bf 1}_j}}f(z)
\ket{{\bu_{\ba +{\bf 1}_i + {\bf 1}_j + {\bf 1}_n}}}
\bra{{\bu_{\ba +{\bf 1}_i + {\bf 1}_j + {\bf 1}_n}}}E\ket{\bu_\ba}
=\sum_{n=1}^N \bra{\bu_{\ba +{\bf 1}_i+{\bf 1}_j}} E \ket{\bu_{\ba - {\bf 1}_n}}
\bra{\bu_{\ba-{\bf 1}_n}}f(z)\ket{\bu_\ba}
$$
For generic $u_1,\ldots,u_N$, it is easy to see that
$\bra{{\bu_{\ba +{\bf 1}_i + {\bf 1}_j + {\bf 1}_n}}}E\ket{\bu_\ba}=0$ comparing
the coefficients of the delta functions. As far as the actions of $e_m$ are well-defined, the Serre relations $E=0$
is valid in the limiting case, too.
\end{proof}

The following lemma is dual to Lemma \ref{sub}.
\begin{lem}
Let $A\subset \Z^N$ be a subset such that
\begin{itemize}
\item for all $\ba\in A$, $\ba'\in\Z^N$ the matrix coefficients $\bra{\bu_{\ba'}} e(z) \ket{\bu_\ba}$ and
$\bra{\bu_{\ba'}} f(z) \ket{\bu_\ba}$ are well-defined,
\item for all $\ba\notin A$, ${\bf b}\in A$ the matrix coefficients
$\bra{\bu_{\bf b}} e(z) \ket{\bu_\ba}$ and $\bra{\bu_{\bf b}} f(z) \ket{\bu_\ba}$ vanish.
\end{itemize}
Then formulas \eqref{ereg}, \eqref{freg} and \eqref{psireg} define a structure
of $\hE$-module on $\mathrm{span}\{\bu_\ba\}_{\ba\in A}$.
\end{lem}
\begin{proof}
Similar to the proof of Lemma \ref{sub}.
\end{proof}

In the following lemma we check that for generic values of parameters
$u_1,\dots,u_N$ the tensor product $V(u_1)\T\dots\T V(u_N)$ is well-defined.
\begin{lem}
Let $u_1,\dots,u_N\in\C$ be
some numbers with the property
\begin{equation}\label{cond}
\frac{u_i}{u_j}\ne q_1^k \text{ for all } 1\le i< j\le N, k\in\Z.
\end{equation}
Then the comultiplication rule \eqref{tre}, \eqref{trf}, \eqref{trpsi}
define the structure of
$\hE$-module on the tensor product $V(u_1)\T \dots\T V(u_N)$.
\end{lem}
\begin{proof}
The action of $\psi^\pm(z)$ is obviously well-defined. They have only simple poles
because of the condition \eqref{cond}.
We check
that the action of $e(z)$ is also well-defined (the case of $f(z)$ is
similar).

By definition we have
\begin{multline}\label{epsi}
e(z) ([u_1]_{i_1}\T\dots\T [u_N]_{i_N})\\
=\sum_{s=1}^N
\left(\prod_{l=1}^{s-1} \gamma_{i_l,u_l}(q_1^{i_s}u_s)\right)
\delta(q_1^{i_s}u_s/z)
[u_1]_{i_1}\T\dots \T [u_{s-1}]_{i_{s-1}}\T  [u_s]_{i_s+1}\T [u_{s+1}]_{i_{s+1}}\T \dots\T  [u_N]_{i_N},
\end{multline}
After substitution, the denominators take the form
$1-q_1^ku_t/u_s$, $k\in\Z$, which do not vanish because of the condition \eqref{cond}.
\end{proof}

\subsection{Submodules of tensor products}
We now consider the tensor product of modules $V(u)$, where the evaluation
parameters form a geometric progression with ratio $q_2$.

Let $V^N(u)$ be the $\hE$-module defined by
\begin{align*}
V^N(u)=V(u)\otimes V(uq_2^{-1})\otimes\cdots\otimes V(uq_2^{-N+1}).
\end{align*}
Set
\begin{align*}
\mathcal P^N=\{\la=(\la_1,\ldots,\la_N)\in\mathbb Z^N|\ \la_1\geq\cdots\geq\la_N\}.
\end{align*}
Let $W^N(u)\hk V^N(u)$ be the subspace spanned by the vectors
\begin{align}\label{ketla}
&\ket{\la}_u=[u]_{\la_1}\otimes [uq_2^{-1}]_{-1+\la_2}
\otimes\cdots\otimes [uq_2^{-N+1}]_{-N+1+\la_N}
\end{align}
where $\la\in\mathcal P^N$. In what follows if the value of $u$ is clear from the context
we abbreviate $\ket{\la}_u=\ket{\la}$.

\begin{lem}
$W^N(u)$ is a level $(1,1)$ submodule of $V^N(u)$.
\end{lem}
\begin{proof}
We prove that $W^N(u)$ is invariant with respect to $e_i$. The case of $f_i$ is similar.

Recall the comultiplication rule
\[
\triangle e(z)=e(z)\T 1 + \psi^-(z)\T e(z).
\]
Since $e(z)[u]_j$ is proportional to $[u]_{j+1}$ it suffices to check that
\[
\psi^-(z) [u]_j\T e(z) [uq_2^{-1}]_{j-1}=0.
\]
By definition,
\[
\psi^-(z) [u]_j\T e(z) [uq_2^{-1}]_{j-1}=
\frac{(1-q_1^{-j}q_3^{-1}z/u)(1-q_1^{-j}q_2^{-1}z/u)}{(1-q_1^{-j}z/u)(1-q_1^{-j+1}z/u)}
\delta(q_1^{j-1}q_2^{-1}u/z),
\]
which vanishes because of \eqref{delta}. The lemma is proved.
\end{proof}

In the following proposition we write down the action of generators
of $\hE$ on $\ket{\la}$ explicitly.
We introduce the notation
$$
\la\pm {\bf 1}_j=(\la_1,\cdots,\la_j\pm 1,\cdots,\la_N).
$$
\begin{prop}\label{mcoef}
The action of $e(z)$ is given by the formula
\begin{equation}\label{e(z)}
(1-q_1)e(z) \ket{\la}
=\sum_{i=1}^N
\prod_{j=1}^{i-1}
\frac{(1-q_1^{\la_i-\la_j}q_3^{i-j-1})(1-q_1^{\la_i-\la_j+1}q_3^{i-j+1})}
{(1-q_1^{\la_i-\la_j}q_3^{i-j})(1-q_1^{\la_i-\la_j+1}q_3^{i-j})}
\delta(q_1^{\la_i}q_3^{i-1}u/z)\ket{\la+{\bf 1}_i}.
\end{equation}
The action of $f(z)$ is given by
\begin{equation}\label{f(z)}
-(1-q_1^{-1})f(z) \ket{\la}
=\sum_{i=1}^N
\prod_{j=i+1}^N
\frac{(1-q_1^{\la_j-\la_i+1}q_3^{j-i+1})(1-q_1^{\la_j-\la_i}q_3^{j-i-1})}
{(1-q_1^{\la_j-\la_i+1}q_3^{j-i})(1-q_1^{\la_j-\la_i}q_3^{j-i})}
\delta(q_1^{\la_i-1}q_3^{i-1}u/z)\ket{\la-{\bf 1}_i}.
\end{equation}
For $\psi$ operators one has
\begin{gather}\label{psi+}
\psi^+(z)\ket{\la}
=\prod_{i=1}^N
\frac{(1-q_1^{\la_i}q_3^iu/z)(1-q_1^{\la_i-1}q_3^{i-2}u/z)}{(1-q_1^{\la_i}q_3^{i-1}u/z)(1-q_1^{\la_i-1}q_3^{i-1}u/z)} \ket{\la},\\
\label{psi-}
\psi^-(z)\ket{\la}
=\prod_{i=1}^N
\frac{(1-q_1^{-\la_i}q_3^{-i}z/u)(1-q_1^{-\la_i+1}q_3^{-i+2}z/u)}
{(1-q_1^{-\la_i}q_3^{-i+1}z/u)(1-q_1^{-\la_i+1}q_3^{-i+1}z/u)}
\ket{\la}.
\end{gather}
\end{prop}
\begin{proof}
Follows from the comultiplication rules and the definition of the modules $V(u)$.
\end{proof}

\section{Semi-infinite construction}
In this section we construct the Fock modules $\F(u)$ by using the inductive limit of
certain subspaces in the finite tensor products of vector representations.
In order to construct the inductive limit consistently we need to modify the operators 
$f(z)$ and $\psi^\pm(z)$, and
this modification results in the nontrivial level $(1,q_2)$ in the inductive limit.
\subsection{Modified operators}
Recall that for any $N\ge 1$ the basis of the space $W^N(u)$ is labeled by the sequences
$\la=(\lambda_1,\dots ,\lambda_N)\in\mathcal P^N$.
The corresponding vectors $\ket{\lambda}\in W^N(u)$ are given by formula \eqref{ketla}.
Set
\begin{align*}
\mathcal P^{N,+}=\{\la\in\mathcal P^N|\la_N\geq0\},
\end{align*}
and define
$W^{N,+}(u)$ to be the subspace of $W^N(u)$ spanned by the vectors $\ket\la$ for $\la\in\mathcal P^{N,+}$.
Our goal in this section is to construct a semi-infinite tensor product
$\F(u)$ as an inductive limit of $W^{N,+}(u)$, and to endow it with a structure of $\hE$-module.
For this purpose we need to vary $N$, and take the limit $N\rightarrow\infty$.
Our strategy is as follows: let $\tau_N:\mathcal P^{N,+}\rightarrow\mathcal P^{N+1,+}$
be the mapping given by
\begin{align*}
\tau_N(\la)=(\lambda_1,\dots ,\lambda_N,0),
\end{align*}
and induce the embedding $\tau_N:W^{N,+}(u)\hk W^{N+1,+}(u)$. Define the inductive limit
\begin{align}
\label{indlim}
\F(u)=\lim_{N\rightarrow\infty}W^{N,+}(u).
\end{align}
This space is spanned by the vectors $\ket{\la}\ (\la\in\mathcal P^+)$ where
the sets of infinite partitions $\mathcal P,\mathcal P^+$ are defined by
\begin{align*}
&\mathcal P=\{\la=(\la_1,\la_2,\ldots)|\la_i\geq\la_{i+1},\la_i\in\mathbb Z\},\\
&\mathcal P^+=\{\la\in\mathcal P|\la_i=0\ \hbox{\rm for sufficiently large $i$}\}.
\end{align*}
In what follows we refer to $\F(u)$ as Fock module.
We define an action of $\hE$ on $\F(u)$ in the following way.
Consider the  operators acting
on the space $W^N(u)$:
\begin{gather}\label{efN}
e^{[N]}(z)=e(z),\quad f^{[N]}(z)=\frac{1-q_2q_3^Nu/z}{1-q_3^Nu/z} f(z),\\
\label{+-N}
\psi^{+[N]}(z)=\frac{1-q_2q_3^Nu/z}{1-q_3^Nu/z}\psi^+(z),\
\psi^{-[N]}(z)=q_2\frac{1-q_2^{-1}q_3^{-N}z/u}{1-q_3^{-N}z/u}\psi^-(z).
\end{gather}

\medskip
\begin{rem}
The action of $f(z)$ splits into a sum of delta functions. By definition
the change caused by the multiplication of the rational function
\[
\beta_N(z)=\frac{1-q_2q_3^Nu/z}{1-q_3^Nu/z}
\]
is that each delta function is multiplied by the value of $\beta_N(z)$
at its support. The changes in $\psi^+(z)$ and $\psi^-(z)$
are the multiplication by $\beta_N(z)$ as a series in $z^{-1}$ and
that in $z$ respectively. Since $\beta_N(z)$ has no pole at $z=\infty$ or $z=0$,
the regularity of the series is not violated in both cases.

Another point is why we choose $\beta_N(z)$ to multiply. This factor has the following meaning.
Consider the eigenvalue of $\psi^+(z)$ on the tensor component $[q_2^{-N+1}u]_{\la_N-N+1}$
when $\la_N=0$. It has $4$ factors, say of the form $\frac{\al^{(1)}_N(z)\al^{(2)}_N(z)}{\al^{(3)}_N(z)\al^{(4)}_N(z)}$.
In fact, we have
\begin{align*}
\beta_N(z)=\frac{\al^{(4)}_N(z)}{\al^{(2)}_N(z)}=\frac{\al^{(1)}_{N+1}(z)}{\al^{(3)}_{N+1}(z)}.
\end{align*}
Namely, we have removed a part of the  factors from the tail
(we say the $N$-th tensor component is in the tail if $\la_N=0$); those which will disappear
when we extend the tail from $N$ to $N+1$.
\end{rem}
\medskip

It turns out that the operators $x^{[N]}(z)$, $x=e,f,\psi^\pm$
are stable and define an $\hE$-module structure on $\F(u)$.
Let us give precise definitions.


First we prepare
\begin{lem}\label{IND}
Suppose that for $\la\in\mathcal P^{N,+}$ the equality $\la_N=0$ is valid.
Then, for $x=e,f,\psi^+,\psi^-$ we have
$x^{[N]}(z)\ket{\la}\in W^{N,+}(u)$ and
\[
\tau_N\left(x^{[N]}(z)\ket{\la}\right)=
x^{[N+1]}(z)\tau_N\left(\ket{\la}\right).
\]
\end{lem}

\begin{proof}
For $x=e$ our lemma is trivial
because in the right hand side the $(N+1)$-st term in the comultiplication of $e(z)$
acts trivially. Let $x=\psi^+$.
We need to prove that the eigenvalue of the operator
\[
\beta_N(z)\psi^+(z) \text{ on the vector } \ket{\la_1,\dots,\la_N}
\]
coincides with that of the operator
\[
\beta_{N+1}(z)\psi^+(z) \text{ on the vector } \ket{\la_1,\dots,\la_N,0}.
\]
Since
$$
\psi^+(z)[uq_2^{-N}]_{-N}=\frac{(1-q_1^{-N}q_3q_2^{-N}u/z)(1-q_1^{-N}q_2^{-N+1}u/z)}
{(1-q_1^{-N}q_2^{-N}u/z)(1-q_1^{-N-1}q_2^{-N}u/z)}[uq_2^{-N}]_{-N}
$$
the statement follows from the equality
\[
\beta_N(z)=\beta_{N+1}(z)\frac{(1-q_3^{N+1}u/z)(1-q_2q_3^Nu/z)}
{(1-q_3^Nu/z)(1-q_3^{N+1}q_2u/z)}.
\]
The case $x=\psi^-$ is similar.
For $x=f$ the comparison of the $i$-th terms in the left and right hand sides is similar
for $1\leq i\leq N-1$. It is easy to see that the rest of the terms, i.e., $i=N$ in the left hand side
and $i=N,N+1$ in the right hand side, are zero. So, the equality for $x=f$ follows.
\end{proof}

\subsection{Fock modules}
We now endow each space $\F(u)$ with a structure of $\hE$-module. For any
$\la=(\la_1,\la_2,\ldots)\in\mathcal P^+$ we set
\begin{equation}\label{infty}
x(z)\ket{\la}=\lim_{N\rightarrow\infty}x^{[N]}(z)\ket{\la_1,\ldots,\la_N},
\end{equation}
where $x=e,f,\psi^+,\psi^-$ and the right hand side is
considered as an element of $\F(u)$ via \eqref{indlim}.

\begin{thm}\label{F(u)}
Formula \eqref{infty} endows $\F(u)$ with the structure of level $(1,q_2)$ $\hE$-module.
\end{thm}
\begin{proof}
We have to check that all relations of $\hE$ are satisfied.
The only non-trivial check is the commutation relation
\[
[e(z),f(w)]=\frac{\delta(z/w)}{g(1,1)}(\psi^+(z)-\psi^-(z)).
\]
We prove this relation using the structure of $\hE$-module on $W^N(u)$.
For given $\ket{\la}\in \F(u)$ choose $N$ large enough so that $\la_{N-1}=\la_N=0$.
Set $\la^{[N]}=(\la_1,\ldots,\la_N)$.
Because of Lemma \ref{IND}, it suffices to prove the equality
\begin{align}\label{EF}
[e^{[N]}(z),f^{[N]}(w)]\ket{\la^{[N]}}
=\frac{\delta(z/w)}{g(1,1)}(\psi^{+[N]}(z)-\psi^{-[N]}(z))\ket{\la^{[N]}}.
\end{align}
We start with the equality in $W^N(u)$:
\[
[e(z),f(w)]\ket{\la^{[N]}}
=\frac{\delta(z/w)}{g(1,1)}(\psi^+(z)-\psi^-(z))\ket{\la^{[N]}}.
\]
This is an equality for the coefficients of $\delta(z/w)\ket{\la^{[N]}}$.
The coefficients are Laurent series in $z$.
In the left hand side the coefficient is a sum of delta functions:
\[
\hbox{\rm LHS}=\sum_{i=1}^Nc_i\delta(q_1^{\la_i}q_3^{i-1}u/z)
+\sum_{i=1}^Nc'_i\delta(q_1^{\la_i-1}q_3^{i-1}u/z).
\]
The right hand side is expressed in terms of a rational function $a_N(z)$:
the operator $\psi^+(z)$ has an eigenvalue on $\ket{\la^{[N]}}$ which is equal to
the series expansion of $a_N(z)$ in $z^{-1}$, while the eigenvalue of $\psi^-(z)$ is
the series expansion of the same rational function $a_N(z)$ in $z^{-1}$.
The rational function $a_N(z)$ is regular at both $z=0$ and $z=\infty$.
The equality implies that $a_N(z)$ has the only simple poles at $z=q_1^{\la_i}q_3^{i-1}u$
with the residue $c_i$, and at $z=q_1^{\la_i-1}q_3^{i-1}u$ with the residue $c'_i$.
Moreover, since $\la_N=0$, it has a zero at $z=q_3^Nu$.

Our aim is to prove \eqref{EF}. To obtain the right hand side of this equality,
we expand $a_N(z)\beta_N(z)$ in $z^{-1}$ for $\psi^+(z)$, and in $z$ for $\psi^-(z)$,
and then take the difference.
We get a Laurent series as the difference of these two expansions.
Note that the only poles of $a_N(z)\beta_N(z)$ in $\mathbb C\sqcup\{\infty\}$
are still at $z=q_1^{\la_i}q_3^{i-1}u$ and $z=q_1^{\la_i-1}q_3^{i-1}u$ after the multiplication of $\beta_N(z)$.
Therefore, the $n$-th Fourier coefficient of this series is calculated by taking the sum of the residues
of $a_N(z)\beta_N(z)z^{-n-1}dz$ at these poles.
On the other hand the same procedure applied to the partial fraction of $a_N(z)\beta_N(z)$
gives rise to the change of the coefficients of the delta functions by
multiplication of the values of $\beta_N(z)$ at their supports.
Therefore we have the equality of two Laurent series.

We give a remark on the proof of the Serre relations. For $E=[e_0,[e_1,e_{-1}]]$ the proof follows from the
case of the finite tensor product. For $F=[f_0,[f_1,f_{-1}]]$, it is the same because the modification of the
actions is such that $\bra{\la-{\bf 1}_i-{\bf 1}_j-{\bf 1}_n}F\ket\la$ changes by a constant multiple.
\end{proof}

\begin{cor}
\label{Fu}
The non-zero matrix coefficients of the action of the generators on $\F(u)$
are given by the formulas:\\
For $e(z):$
$$
(1-q_1)\bra{\la+{\bf 1}_i}e(z) \ket{\la}
=\prod_{j=1}^{i-1}
\frac{(1-q_1^{\la_i-\la_j}q_3^{i-j-1})(1-q_1^{\la_i-\la_j+1}q_3^{i-j+1})}
{(1-q_1^{\la_i-\la_j}q_3^{i-j})(1-q_1^{\la_i-\la_j+1}q_3^{i-j})}
\delta(q_1^{\la_i}q_3^{i-1}u/z).
$$
For $f(z):$
\begin{multline*}
-(1-q_1^{-1})\bra{\la - {\bf 1}_i}f(z) \ket{\la}\\
=\frac{1-q_1^{\la_{i+1}-\la_i}}{1-q_1^{\la_{i+1}-\la_i+1}q_3}
\prod_{j=i+1}^\infty
\frac{(1-q_1^{\la_j-\la_i+1}q_3^{j-i+1})((1-q_1^{\la_{j+1}-\la_i}q_3^{j-i})}
{(1-q_1^{\la_{j+1}-\la_i+1}q_3^{j-i+1})(1-q_1^{\la_j-\la_i}q_3^{j-i})}
\delta(q_1^{\la_i-1}q_3^{i-1}u/z).
\end{multline*}
For $\psi^\pm(z):$
\begin{gather*}
\psi^+(z)\ket{\la}
=\frac{1-q_1^{\la_1-1}q_3^{-1}u/z}{1-q_1^{\la_1}u/z}
\prod_{i=1}^\infty
\frac{(1-q_1^{\la_i}q_3^iu/z)(1-q_1^{\la_{i+1}-1}q_3^{i-1}u/z)}{(1-q_1^{\la_{i+1}}q_3^iu/z)(1-q_1^{\la_i-1}q_3^{i-1}u/z)}
\ket{\la},\\
\psi^-(z) \ket{\la}
=q_2\frac{1-q_1^{-\la_1+1}q_3z/u}{1-q_1^{-\la_1}z/u}
\prod_{i=1}^\infty
\frac{(1-q_1^{-\la_i}q_3^{-i}z/u)(1-q_1^{-\la_{i+1}+1}q_3^{-i+1}z/u)}{(1-q_1^{-\la_{i+1}}q_3^{-i}z/u)(1-q_1^{-\la_i+1}q_3^{-i+1}z/u)}
\ket{\la}.
\end{gather*}
\end{cor}


We define the vacuum vector $\ket{\la^0}$ with $\la^0_i=0$.
Let $\psi^\pm_{\emptyset}(z)$ be the expansions
of a rational function $\frac{1-q_2z}{1-z}$
as a series in $z^{\pm 1}$. Explicitly,
$$
\psi^+_{\emptyset}(u/z)=\frac{1-q_2z}{1-z},\quad \psi^-_{\emptyset}(u/z)=q_2\frac{1-q_2^{-1}z^{-1}}{1-z^{-1}}.
$$
Then $\psi^\pm_{\emptyset}(u/z)$ are
the eigenvalues of $\psi^\pm(z)$ on $\ket{\la^0}$, i.e.
\[
\psi^\pm(z)\ket{\la^0}=\psi^\pm_{\emptyset}(u/z)\ket{\la^0}.
\]

For a partition $\la$ we denote by $\bra{\la}\psi^\pm(z)_i\ket{\la}$ the eigenvalue of
the series $\psi^\pm(z)$ on the vector $[uq_2^{-i+1}]_{\la_i-i+1}\in V(uq_2^{-i+1})$, i.e.
$$
\psi^\pm(z)[uq_2^{-i+1}]_{\la_i-i+1}=\bra{\la}\psi^\pm(z)_i\ket{\la}[uq_2^{-i+1}]_{\la_i-i+1}.
$$
The index $i$ in $\psi^\pm(z)_i$ indicates the component $V(uq_2^{-i+1})$, and the shift
$\la_i\mapsto\la_i-i+1$ as well.
Similarly, we introduce the matrix coefficients
$\bra{\la+{\bf 1}_i}e(z)_i\ket{\la}$ and $\bra{\la}f(z)_i\ket{\la+{\bf 1}_i}$.
Then the formulas from Corollary \ref{Fu} can be rewritten in the following way:
\begin{gather*}
\bra{\la}\psi^\pm(z)\ket{\la}=\psi^\pm_{\emptyset}(u/z)\prod_{i\ge 1}
\frac{\bra{\la}\psi^\pm(z)_i\ket{\la}}{\bra{\la^0}\psi^\pm(z)_i\ket{\la^0}},\\
\bra{\la+{\bf 1}_i}e(z)\ket{\la} =\bra{\la+{\bf 1}_i}e(z)_i\ket{\la}\prod_{j=1}^{i-1} \bra{\la}\psi^-(z)_j\ket{\la},\\
\bra{\la}f(z)\ket{\la +{\bf 1}_i} = \bra{\la}f(z)_i\ket{\la+{\bf 1}_i}\times\psi^+_\emptyset (q_3^i u/z)
\prod_{j=i+1}^\infty \frac{\bra{\la}\psi^+(z)_j\ket{\la}}{\bra{\la^0}\psi^+(z)_j\ket{\la^0}}.
\end{gather*}

\begin{cor}
The module $\F(1)$ is isomorphic to the module constructed in \cite{FT}, Theorem $3.5$.
\end{cor}
\begin{proof}
We recall that in \cite{FT} the operators $e_i$, $f_j$ and $\psi^\pm_i$  were
constructed on the space $\F$ with basis $[\la]$ labeled by infinite partitions $\la$.
This space is defined as the direct sum of localized equivariant $K$ groups of Hilbert schemes
of points of $\C^2$. The matrix coefficients of these operators are given in
Proposition 3.7, \cite{FT}.
We prove that $\F(1)\simeq \F$.
To this end we identify $t_1=q_1$, $t_2=q_3$
and do a change of basis as follows. Consider  constants $c_\la$ defined by
$c_{\la^0}=1$ and
\begin{equation}\label{cla}
\frac{c_{\la+{\bf 1}_i}}{c_\la}=(1-q_1q_3)
\prod_{j=i}^\infty
\frac{1-q_1^{\la_{j+1}-\la_i}q_3^{j-i+1}}{1-q_1^{\la_j-\la_i}q_3^{j-i+1}}
\prod_{j=1}^{i-1}
\frac{1-q_1^{\la_j-\la_i-1}q_3^{j-i-1}}{1-q_1^{\la_j-\la_i-1}q_3^{j-i}}.
\end{equation}
It is straightforward to check that $c_\la$ are well-defined, i.e. that
the right hand sides $d_{\la,i}$ of \eqref{cla} satisfy
\[
d_{\la+{\bf 1}_i,k}d_{\la,i}=d_{\la+{\bf 1}_k,i}d_{\la,k}.
\]
Another straightforward check shows that the linear map $\F\to \F(1)$, $[\la]\mapsto c_\la\ket{\la}$
is the isomorphism of modules of $\hE$.
\end{proof}

\medskip
\begin{rem}
Strictly speaking, in \cite{FT} the authors proved that the operators
$e(z)$, $f(z)$ and $\psi^\pm(z)$ acting on $\F$ satisfy the relations of the Ding-Iohara algebra
$\hE'$.
However, since Serre relations are satisfied on $\F(1)$, the representation
$\F$ factors through the surjection $\hE'\to\hE$.
\end{rem}
\medskip

\section{Macdonald polynomials and spherical DAHA}
In this section we establish a link between the spherical double affine Hecke 
algebra and quantum continuous $\gli$.
Throughout the section we consider the algebra $\hE$ over the field of rational functions
$\C(q_1,q_3)$.
\subsection{Macdonald polynomials}
Our basic reference in this section is the Macdonald's book \cite{M}.
However we use the Laurent polynomials version of Macdonald polynomials.

The Macdonald operators $D^r_N$ are mutually commuting $q$-difference operators
acting on the ring of symmetric Laurent polynomials
$\C(q,t)[x_1^{\pm 1} ,\dots,x_N^{\pm 1}]^{\mathfrak{S}_N}$,
where $\mathfrak{S}_N$ denotes the symmetric group of $N$ letters. These operators are given by
the formula
$$
D^r_N=\sum_{|I|=r} A_I(x;t)T_I,
$$
where $I\subset\{1,\dots,N\}$ runs over subsets of cardinality $r$,
\begin{gather*}
A_I(x;t)=t^{r(r-1)/2}\prod_{i\in I, j\notin I} \frac{tx_i-x_j}{x_i-x_j},\\
T_I=\prod_{i\in I} T_{q,x_i}
\end{gather*}
and $(T_{q,x_i}f)(x_1,\dots,x_N)=f(x_1,\dots,qx_i,\dots,x_N)$.
Let $D_N(X;q,t)=\sum_{r=0}^N D_N^rX^r$ be their generating function.

The Macdonald polynomials $P_\la$ form a basis in the space
of symmetric polynomials. They are uniquely characterized by the following defining properties
\begin{gather*}
D_N(X;q,t)P_\la=\prod_{i=1}^N (1+Xq^{\la_i}t^{N-i})\cdot P_\la,\\
P_\la=m_\la+\sum_{\mu<\la}u_{\la\mu}m_\mu \qquad u_{\la\mu}\in\C(q,t),
\end{gather*}
where $m_\la$ denotes the monomial symmetric function and we write $\mu < \la$ if
$\mu\ne\la$ and $\mu_1+\dots +\mu_i\le \la_1+\dots +\la_i$ for $i=1,\dots,N$.

We define $P_\la(x;q,t)$ for partitions $\la\in\mathcal P^N$ with possibly negative entries by the formula:
\begin{align*}
P_{\la}(x;q,t)
=\prod_{i=1}^Nx_i^{\la_N}
\cdot
P_{\la_1-\la_N,\cdots,\la_{N-1}-\la_N,0}(x;q,t)\,.
\end{align*}
In what follows the following Macdonald operators will be of special importance for us:
\begin{align*}
&D_N^1(q,t)
=\sum_{i=1}^N\prod_{j(\neq i)}\frac{tx_i-x_j}{x_i-x_j}
T_{q,x_i}\,,
\\
&D_N^{-1}(q,t)=D_N^1(q^{-1},t^{-1}).
\end{align*}

Set $W^N=W^N(1)$ (i.e. $u=1$).
\begin{prop}\label{Macdonald}
Choose $q_1=q,q_2=q^{-1}t,q_3=t^{-1}$.
Under the isomorphism of vector spaces
\begin{align*}
W^N\overset{\sim}
{\longrightarrow} \C[x_1^{\pm1},\cdots,x_N^{\pm1}]^{\mathfrak{S}_N},
\quad
\ket{\la}\mapsto P_\la(x)\,,
\end{align*}
we have the identification
\begin{align*}
&(1-q_1)e_0=\text{multiplication by $\sum_{i=1}^Nx_i$},
\\
&-(1-q_1^{-1})f_0=\text{multiplication by
$\sum_{i=1}^Nx_i^{-1}$},
\\
&q_3^{N-1}(1-q_2)(1-q_3)\psi^+_1=D^1_N(q,t),
\\
&
q_3^{-N+1}(1-q_2^{-1})(1-q_3^{-1})\psi^-_{-1}=D^{-1}_N(q,t)\,.
\end{align*}
\end{prop}
\begin{proof}
We first look at the operator $\psi^+_1$ (the $\psi^-_{-1}$-case is  similar).
Formula \eqref{psi+} gives
\[
\psi^+_1\ket{\la}=q_3^{N-1}(1-q_2)(1-q_3)\sum_{i=1}^N q_1^{\la_i} q_3^{i-N}\ket{\la},
\]
which agrees with the formula for the eigenvalues of the operator $D^1_N$ (see \cite{M}).

In order to prove that $(1-q_1)e_0$ acts as $\sum x_i$ we compare the matrix
coefficients of $(1-q_1)e_0$ in the basis $\ket{\la}$ and of $\sum x_i$ in
the basis $P_\la$. The latter is given by the Pieri formulas.
From \eqref{e(z)} we obtain
\begin{equation}\label{m}
(1-q_1)e_0\cdot \ket{\la}
=\sum_{i=1}^N
\prod_{j=1}^{i-1}
\frac{(1-q_1^{\la_i-\la_j}q_3^{i-j-1})(1-q_1^{\la_i-\la_j+1}q_3^{i-j+1})}
{(1-q_1^{\la_i-\la_j}q_3^{i-j})(1-q_1^{\la_i-\la_j+1}q_3^{i-j})}
\cdot\ket{\la+{\bf 1}_i}.
\end{equation}
The matrix coefficients as above coincide with the Pieri rule
formulas (see \cite{M}). Similarly one proves that $-(1-q_1^{-1})f_0$ acts as multiplication
by $\sum_{i=1}^N x_i^{-1}$.
\end{proof}

\subsection{Spherical double affine Hecke algebras}
We first recall the definition of the double affine Hecke algebra of type
$GL_N$ (DAHA for short) \cite{C}. This is a $\C(q,v)$-algebra generated by elements
$T^{\pm 1}_i$, $X^{\pm 1}_j$ and $Y^{\pm 1}_j$ for $1\le i\le N-1, 1\le j\le N$,
subject to the following relations
\begin{gather*}
(T_i+v^{-1})(T_i-v)=0,\quad T_i T_{i+1} T_i=T_{i+1} T_i T_{i+1},\\
T_i T_k=T_k T_i \text{ if } |i-k|>1,\\
X_j X_k=X_k X_j, \quad Y_j Y_k=Y_k Y_j,\\
T_i X_i T_i=X_{i+1}, \quad T^{-1}_i Y_i T^{-1}_i=Y_{i+1},\\
T_i X_k=X_k T_i,\quad T_i Y_k=Y_k T_i \text { if } k\ne i,i+1,\\
Y_1 X_1 \dots X_N=q X_1 \dots X_NY_1,\\
X^{-1}_1 Y_2=Y_2 X^{-1}_1 T^{-2}_1.
\end{gather*}
Denote this algebra by $\ddot H_N$. Let $S\in \ddot H_N$ be the idempotent given by
$$
S=\frac{1}{[N]!} \sum_{w\in \mathfrak{S}_N}v^{l(w)} T_w,\quad T_w=T_{i_1} \dots T_{i_r},
$$
for a reduced decomposition $w=s_{i_1} \dots s_{i_r}$.
Here $s_i$ denotes the transposition $(i,i+1)$, $l(w)$ is the length of $w$
and
$$[N]!=\prod^N_{i=1}[i],\quad [i]=\frac{v^{2i}-1}{v^2-1}.$$
The algebra $S\ddot H_N S$ is called the spherical DAHA and is denoted by $S\ddot H_N$.
In \cite{SV1}
Schiffmann and Vasserot defined the elements $P^N_{0,m}, P^N_{m,0}\in S\ddot H_N$ ($m\in\Z$)
by the formulas:
\begin{gather*}
P^N_{0,l}=S\sum^N_{i=1} Y^l_iS,\quad P^N_{0,-l}=q^l S \sum^N_{i=1}Y^{-l}_i S,\\
P^N_{l,0}=q^l S \sum_{i=1}^N X^l_i S,\quad P^N_{-l,0}=S \sum_{i=1}^N X^{-l}_i S
\end{gather*}
with $l>0$.
They proved that these elements generate $S\ddot H_N$.
We need the following modification of their result:

\begin{lem}\label{4SH}
Four elements $P^N_{0,1},\ P^N_{0,-1},\ P^N_{1,0} \text{ and} \ P^N_{-1,0}$ generate $S\ddot H_N$.
\end{lem}
\begin{proof}
Consider the degeneration $v=1$ of $\ddot H_N$.
(In fact, one has to be careful with such a degeneration since $\ddot H_N$
is defined over the field of rational functions $\C(q,v)$, which may have a pole
at $v=1$. In order to make everything precise, one has to pass to the analogue
of $\ddot H_N$, defined over $\C[q^{\pm 1},v^{\pm 1}]$. We omit the details here
and refer the reader to \cite{SV1}, Section $2$.)
We prove that for  $v=1$ the corresponding spherical DAHA is generated by our $4$ elements.
This  would imply the lemma.
If $v=1$, the idempotent $S$ commutes with $\sum^N_{i=1} X^{\pm 1}_i$ and $\sum^N_{i=1}Y^{\pm 1}_i$
and $X_i Y_j=q^{-\delta_{i,j}}Y_j X_i$.
Let
$$P_{\pm 1}=S\sum_{i=1}^N X^{\pm 1}_i S, \quad Q_{\pm 1}=S\sum_{i=1}^N Y^{\pm 1}_i S.$$
Then we have
\begin{align*}
(ad P_{\pm 1})^k Q_1&=S\sum_{i=1}^N (ad X^{\pm 1}_i)^k Y_i S\\
&=(1-q^{\pm 1})^k S \sum_{i=1}^N X^{\pm k}_i Y_i S
\end{align*}
and for an arbitrary $m\in\Z\setminus\{ 0\}$
\begin{align*}
(ad Q_{\pm 1})^m S\sum_{i=1}^N X^k_i Y_i S&=S\sum_{i=1}^N (ad Y^{\pm 1}_i)^m X^k_i Y_i S\\
&=(1-q^{\pm 1})^{mk} S\sum_{i=1}^N X^k_i Y^{1\pm m}_i S.
\end{align*}
We thus obtain that all $P^N_{m,0}$ and $P^N_{0,m}$ can be obtained as linear combinations of products of
$P^N_{\pm 1,0}$ and $P^N_{0,\pm 1}$.
\end{proof}

Similar result holds for the quantum continuous $\gli$.

\begin{lem}\label{4DI}
For $c^{\pm}\in \C^\times$, the algebra  $\hE/\langle\psi^{\pm}_0-c^{\pm}\rangle $
is generated
by four elements $e_0,\ \psi^+_1,\ f_0,\ \psi^-_1$.
\end{lem}

\begin{proof}
Follows directly from relations in $\hE$. For example, the $z^{\mp 1} w^{-m}$
terms of the relations  \eqref{rel2}  give
$$[\psi^+_1,e_m]=c^+\left(\sum_{j=1}^3 q_j-\sum_{j=1}^3
q^{-1}_j\right)e_{m+1},$$
$$[\psi^-_{-1},e_m]=c^-\left(\sum_{j=1}^3 q_j-\sum_{j=1}^3 q^{-1}_j\right)e_{m-1}.$$
Therefore, all $e_m,\ m\in\Z$ can be obtained via $e_0$, $\psi^+_1$ and $\psi^-_{-1}$.
Similarly one gets $f_m,\ m\in\Z, \ \psi^+_n,\ \psi^-_{-n},\ n\ge 0$.
\end{proof}

\begin{thm}\label{DISH}
For any $N$ there exists a surjective homomorphism of algebras
$\hE\to S\ddot H_N$,
where the parameters $q,v$ of $S\ddot H_N$ are related to $q_1,q_3$ by $q=q_1$,
$v^2=q^{-1}_3$.
\end{thm}

\begin{proof}
We recall that the algebra $S\ddot H_N$ can be faithfully represented on the space
$$\C(q,t)[x^{\pm 1}_1,\dots ,x^{\pm 1}_N]^{S_N}$$
with $t=v^2$ in such a way that
$P^N_{\pm 1,0}$ acts as a multiplication by
$\sum^N_{i=1}x^{\pm 1}_i$ and $P^N_{0,\pm 1}$ acts as Macdonald difference operators
$D_N^1$ and $D_N^{-1}$.
Therefore, Proposition \ref{Macdonald} and Lemmas \ref{4SH} and \ref{4DI} show that the assignment
\begin{gather*}
(1-q_1)e_0\mapsto S\sum_{i=1}^N X_i S,\quad -(1-q_1^{-1})f_0\mapsto S\sum_{i=1}^N X^{-1}_i S,\\
q_3^{N-1}(1-q_2)(1-q_3)\psi^+_1\mapsto S\sum_{i=1}^N Y_i S,\quad
q_3^{1-N}(1-q_2^{-1})(1-q_3^{-1})\psi^-_{-1}\mapsto S\sum_{i=1}^N Y^{-1}_1 S
\end{gather*}
extends to the surjective homomorphism of algebras.
\end{proof}

\section{Resonance case}\label{secres}
Let $k\ge 1,r\ge 2$ be positive integers.
In this section we impose the following condition on the parameters $q_1,q_2,q_3$
\begin{equation}\label{rescond}
q_1^{1-r}q_3^{k+1}=1.
\end{equation}
As usual we assume that $q_1q_2q_3=1$. We refer to the condition
\eqref{rescond} as the resonance condition.
We also assume that
$q_1^n q_3^m=1$ if and only if
there exists an integer
$\alpha$ such that $n=(1-r)\al$, $m=(k+1)\al$.

In this section we establish a link between the representations of $\hE$
and ideals in polynomial algebra spanned by the Macdonald polynomials (see
\cite{FJMM1}, \cite{FJMM2}).

\subsection{Finite tensor products}
If the resonance condition holds, the action of $\hE$ on $W^N(u)$ becomes ill-defined
(since the denominators of the formulas determining $e(z)w$ and $f(z)w$
vanish for some vectors $w\in W^N(u)$). We still find a subspace inside $W^N(u)$ on which
the action is well-defined.

Set
\begin{align}\label{fat}
S^{k,r,N}=\{\la\in\mathcal P^N|\la_i-\la_{i+k}\geq r\ (1\leq i\leq N-k)\}.
\end{align}
We call partitions satisfying the condition \eqref{fat},  $(k,r)$-admissible partitions.
Let $W^{k,r,N}(u)\hk W^N(u)$ be the subspace spanned by the vectors $\ket\la$ for $\la\in S^{k,r,N}$.

Induce the actions of
the operators $e(z),f(z),\psi^\pm(z)$ on $W^{k,r,N}(u)$ from those on $W^N(u)$ for the generic values
of the parameters.
\medskip
\begin{rem}
In fact one has to be careful defining the matrix coefficients
$\bra{\la+{\bf 1}_i} e(z) \ket{\la}$ and $\bra{\la -{\bf 1}_i} f(z) \ket{\la}$. Namely
both contain factors of the form $1-q_1^s q_3^l$. If the condition \eqref{rescond}
holds and $s=\al(1-r)$, $l=\al(k+1)$, such a factor vanishes. The prescription is
first to cancel all factors of the form $1-q_1^{1-r}q_3^{k+1}$ (if they appear simultaneously
in the numerator and in the denominator) and then impose the resonance condition.
\end{rem}
\medskip

\begin{lem}
The comultiplication rule makes the subspace $W^{k,r,N}(u)\hk W^N(u)$ into
level $(1,1)$ $\hE$-module.
\end{lem}

\begin{proof}
We need to check that matrix coefficients $\bra{\la+{\bf 1}_i} e(z) \ket{\la}$
($\bra{\la -{\bf 1}_i} f(z) \ket{\la}$) are well-defined provided
$\la + {\bf 1}_i\in S^{k,r,N}$ ($\la - {\bf 1}_i\in S^{k,r,N}$). We check this for $e(z)$
(the case of $f(z)$ is similar). Formula \eqref{e(z)} gives
\begin{equation}
\label{poles}
(1-q_1)\bra{\la+{\bf 1}_i} e(z) \ket{\la}=
\prod_{j=1}^{i-1}
\frac{(1-q_1^{\la_i-\la_j}q_3^{i-j-1})(1-q_1^{\la_i-\la_j+1}q_3^{i-j+1})}
{(1-q_1^{\la_i-\la_j}q_3^{i-j})(1-q_1^{\la_i-\la_j+1}q_3^{i-j})}
\delta(q_1^{\la_i}q_3^{i-1}u/z).
\end{equation}
The denominator  vanishes if for some $m$ satisfying $1\le m\le i-1$
\[
q_1^{\la_i-\la_m}q_3^{i-m}=1 \text{ or } q_1^{\la_i-\la_m+1}q_3^{i-m}=1.
\]
In the first case there exist some positive integer  $\al$ such that
$$
\la_i-\la_m=\al(1-r),\quad i-m=\al(k+1).
$$
This is impossible since $\la + {\bf 1}_i$ is $(k,r)$-admissible. In fact
\[
\la_i + 1\le \la_{i-\al k} - \al r = \la_{m+\al} - \al r \le \la_m - \al r
\]
and hence $\la_i - \la_m\le -1 - \al r < \al (1-r)$.
Now assume $q_1^{\la_i-\la_m+1}q_3^{i-m}=1$. Then there exists a positive integer $\al$ such that
$$
\la_i-\la_m+1=\al(1-r),\quad i-m=\al(k+1).
$$
As above we have $\la_i + 1 - \la_m\le  - \al r < \al (1-r)$, which is a contradiction.


We  show now that the action of $\hE$ preserves the linear span of
vectors corresponding to the $(k,r)$-admissible partitions. In fact, let $\la_i-\la_{i+k}=r$. Then
we need to show that
\begin{equation}\label{braket}
\bra{\la+{\bf 1}_{i+k}} e(z) \ket{\la}=0
\end{equation}
(since $\la+{\bf 1}_{i+k}$ is not $(k,r)$-admissible).
But this zero comes from the factor $1-q_1^{\la_i-\la_j+1}q_3^{i-j+1}$ in the formula \eqref{e(z)}
where $i$ and $j$ are replaced with $i+k$ and $j$ respectively.
The check for $f(z)$ is similar.
\end{proof}

\medskip
\begin{rem}
We recall that in \cite{FJMM2} the vector space spanned by the Macdonald polynomials $P_\la$
with $(k,r)$-admissible partitions $\la$ was considered (see also \cite{Kas}).
It was proved that this space is $S\ddot H_N$
stable. The lemma above gives yet another proof of this statement using the representation
theory of $\hE$.
\end{rem}
\medskip

\subsection{Semi-infinite limit}
In this subsection we define a semi-infinite representation in the resonance case.
The construction is similar to the construction for $\F(u)$,
though certain modifications are needed.

Fix a sequence of integers $\bc=(c_1,\dots,c_{k-1})$ satisfying
$0=c_0\le c_1 \le \cdots\le c_{k-1}\le r$,
and define the tail
\begin{align*}
\la^0_{\nu k+i+1}=-\nu r-c_i\ (\nu\geq0,0\leq i\leq k-1).
\end{align*}
We define the sets of partitions $S^{k,r}_\bc, S^{k,r,N,+}_\bc$ and
the mapping $\tau_{k,r,N}:S_\bc^{k,r,N,+}\to S_\bc^{k,r,N+k,+}$ as follows.
\begin{align*}
&S^{k,r}_\bc=\{\la\in\mathcal P|\ \la_j-\la_{j+k}\ge r\quad (j\ge 1),
\quad\la_j=\la_j^0\ \hbox{\rm for sufficiently large $i$}\},\\
&S^{k,r,N,+}_\bc=\{\la\in S^{k,r,N}|\la_j\ge\la_j^0\ (1\leq j\leq N)\},\\
&\tau_{k,r,N}(\la)=(\la_1,\ldots,\la_N,\la^0_{N+1},\ldots,\la^0_{N+k}).
\end{align*}

Let $W_\bc^{k,r}(u)$ be the space spanned by the vectors $\ket\la$ for $\la\in S^{k,r}_\bc$.
In order to endow it with the structure of $\hE$-module we introduce the subspaces
$W_\bc^{k,r,N,+}(u)\hk W^{k,r,N}(u)$ spanned by the vectors $\ket\la$ for $\la\in S^{k,r,N,+}_\bc$ and
induce the embeddings
\begin{equation}
\tau_{k,r,N}:W_\bc^{k,r,N,+}(u)\to W_\bc^{k,r,N+k,+}(u)
\end{equation}
by the formula $\tau_{k,r,N}\ket{\la}=\ket{\tau_{k,r,N}(\la)}$.
These embeddings give the identification
\begin{equation}\label{rindlim}
W_\bc^{k,r}(u)\simeq \lim_{N\to\infty} W_\bc^{k,r,N+}(u).
\end{equation}

We define an action of $\hE$ on $W_\bc^{k,r}(u)$ in the following way.
Let $\beta_{k,N}(z)$ be the rational function defined by
\[
\beta_{k,N}(z)=\prod_{j=0}^i\frac{1-q_1^{-c_i-\nu-1}q_3^{-\nu+j}u/z}{1-q_1^{-c_i-\nu-1}q_3^{-\nu+j-1}u/z}
\prod_{j=i+1}^{k-1}\frac{1-q_1^{-c_i-\nu}q_3^{-\nu+j+1}u/z}{1-q_1^{-c_i-\nu}q_3^{-\nu+j}u/z}
\]
where $N=\nu k+i+1\ (0\leq i\leq k-1,\nu\geq0)$.
Consider the  operators acting
on the space $W_\bc^{k,r,N,+}(u)$:
\begin{gather}\label{refN}
e^{[k,N]}(z)=e(z),\quad
f^{[k,N]}(z)=\beta_{k,N}(z)f(z),\\
\label{r+-N}
\psi^{+[k,N]}(z)=\beta_{k,N}(z)\psi^+(z),\
\psi^{-[k,N]}(z)=\beta^+_{k,N}(z)\psi^-(z),
\end{gather}
where $\beta^+_{k,N}(z)$ is the expansion of $\beta_{k,N}(z)$ as series in $z$.
These operators turn out to be stable and define the $\hE$-module structure on
$W_\bc^{k,r}(u)$.

In the following, i.e., Lemma \ref{Pkr} and Theorem \ref{Wkr}, the arguments are very similar to those for
Lemma \ref{IND} and Theorem \ref{F(u)}. We omit the proofs for them.
\begin{lem}\label{Pkr}
Suppose that for $\la\in S^{k,r,N,+}_\bc$ the equality
$(\la_{N-k+1},\ldots,\la_N)=(\la^0_{N-k+1},\ldots,\la^0_N)$ is valid.
Then, for $x=e,f,\psi^+,\psi^-$ we have $x^{[k,N]}(z)\ket\la\in W_\bc^{k,r,N,+}(u)$ and
\[
\tau_{k,r,N}\left(x^{[k,N]}(z)\ket\la\right)
= x^{[k,N+k]}(z)\tau_{k,r,N}\left(\ket{\la}\right).
\]
\end{lem}

For $\la\in S^{k,r}_\bc$ we set
\begin{equation}\label{rinfty}
x(z)\ket{\la}=\lim_{N\rightarrow\infty}x^{[k,N]}(z)\ket{\la_1,\ldots,\la_N}.
\end{equation}
where $x=e,f,\psi^+,\psi^-$ and the right hand side is considered as an element
of $W_\bc^{k,r}(u)$ via \eqref{rindlim}.

\begin{thm}\label{Wkr}
Formula \eqref{rinfty} endows $W_\bc^{k,r}(u)$ with the structure of level $(1,q_3^{k})$
$\mc \hE$-module.
\end{thm}

We now write down explicit formulas for the non-zero matrix coefficients
of operators
$e(z)$, $f(z)$ and $\psi^\pm(z)$ acting on $W_\bc^{k,r}(u)$.
For $e(z)$:
$$
(1-q_1)\bra{\la+{\bf 1}_i}e(z) \ket{\la}
=\prod_{j=1}^{i-1}
\frac{(1-q_1^{\la_i-\la_j}q_3^{i-j-1})(1-q_1^{\la_i-\la_j+1}q_3^{i-j+1})}
{(1-q_1^{\la_i-\la_j}q_3^{i-j})(1-q_1^{\la_i-\la_j+1}q_3^{i-j})}
\delta(q_1^{\la_i}q_3^{i-1}u/z).
$$
For $f(z)$:
\begin{multline*}
-(1-q_1^{-1})\bra{\la - {\bf 1}_i}f(z)\ket{\la}\\
=\prod_{j=i+1}^{i+k}\frac{1-q_1^{\la_j-\la_i+1}q_3^{j-i+1}}{1-q_1^{\la_j-\la_i+1}q_3^{j-i}}
\prod_{j=i+1}^\infty
\frac{(1-q_1^{\la_{j+k}-\la_i+1}q_3^{j+k-i+1})(1-q_1^{\la_j-\la_i}q_3^{j-i-1})}
{(1-q_1^{\la_j-\la_i}q_3^{j-i})(1-q_1^{\la_{j+k}-\la_i+1}q_3^{j+k-i})}
\delta(q_1^{\la_i-1}q_3^{i-1}u/z).
\end{multline*}
For $\psi^\pm(z)$:
\begin{gather*}
\psi^+(z)\ket{\la}
=\prod_{i=1}^k \frac{1-q_1^{\la_i}q_3^{i}u/z}{1-q_1^{\la_i}q_3^{i-1}u/z}
\prod_{i=1}^\infty
\frac{(1-q_1^{\la_i-1}q_3^{i-2}u/z)(1-q_1^{\la_{i+k}}q_3^{i+k}u/z)}
{(1-q_1^{\la_{i+k}}q_3^{i+k-1}u/z)(1-q_1^{\la_i-1}q_3^{i-1}u/z)}
\ket{\la},\\
\psi^-(z) \ket{\la}
=q_3^k
\prod_{i=1}^k \frac{1-q_1^{-\la_i}q_3^{-i}z/u}{1-q_1^{-\la_i}q_3^{-i+1}z/u}
\prod_{i=1}^\infty
\frac{(1-q_1^{-\la_i+1}q_3^{-i+2}z/u)(1-q_1^{-\la_{i+k}}q_3^{-i-k}z/u)}
{(1-q_1^{-\la_{i+k}}q_3^{-i-k+1}z/u)(1-q_1^{-\la_i+1}q_3^{-i+1}z/u)}
\ket{\la}.
\end{gather*}

Define series $\varphi^\pm_{\emptyset}(u/z)$ in $z^{\mp 1}$
by
$$
\varphi^+_{\emptyset}(u/z)=\frac{1-q_3u/z}{1-u/z},\quad \varphi^-_{\emptyset}(u/z)=q_3\frac{1-q_3^{-1}z/u}{1-z/u}
$$
Then the formulas above can be rewritten in the following way:
\begin{gather*}
\bra{\la}\psi^\pm(z)\ket{\la}=\prod_{i=0}^{k-1} \varphi^\pm_{\emptyset}(q_3^i q_1^{-c_i}u/z)
\prod_{i\ge 1}
\frac{\bra{\la}\psi^\pm(z)_i\ket{\la}}{\bra{\lambda^0}\psi^\pm(z)_i\ket{\lambda^0}},\\
\bra{\la+{\bf 1}_i}e(z)\ket{\la} = \bra{\la + {\bf 1}_i}e(z)_i\ket{\la}
\prod_{j=1}^{i-1} \bra{\la}\psi^-(z)_j\ket{\la},\\
\bra{\la}f(z)\ket{\la +{\bf 1}_i} = \bra{\la}f(z)_i\ket{\la+{\bf 1}_i}
\prod_{j=i+1}^\infty \frac{\bra{\la}\psi^+(z)_j\ket{\la}}
{\bra{\la^0}\psi^+(z)_j\ket{\la^0}}
\prod_{j=0}^{k-1} \varphi^+_{\emptyset}(q_3^{i+j}q_1^{\lambda^0_{i+j+1}}u/z).
\end{gather*}

\section{Quantum continuous $\gli$: further directions}\label{quot}
As we have shown above there exists a surjective homomorphism form the quantum continuous
$\gli$ to the spherical DAHA of type $GL_N$.
It is therefore natural to expect that $\hE$ is isomorphic to
the stable limit $S\ddot H_\infty=\lim_{N\to\infty} S\ddot H_N$ constructed in
\cite{SV1}, \cite{SV2}.
The related statement is that the representation
obtained as the direct sum of the modules $W^N(u)$ is faithful.

There is also a link between $\hE$ and the so-called shuffle algebra ${\bf S}$
(see \cite{FO}, \cite{FT}). The latter is an algebra generated by variables
$\tilde e_i$, $i\in\Z$.
Let $\hE^+$ be the subalgebras of $\hE$ generated by $e_i$, $i\in\Z$.
We  conjecture that there exist an isomorphism $\pi: \hE^+ \to {\bf S}$,
$\pi^+(e_i)=\tilde e_i$.

In the paper we have studied tensor products of the vector representations of $\hE$. We have
also considered limits of these tensor products, thus constructing
Fock modules $\F(u)$. It is a natural question to study the tensor products of the modules $\F(u)$.
It turns out that the structure of these tensor products is very rich.
In particular there exist
submodules inside $\F(u_1)\T\dots\T \F(u_N)$ whose characters  coincide
with those of the minimal representations of the $\W_N$ algebras up to a trivial factor.

We plan to return to all these questions elsewhere.

\section*{Acknowledgements}
Research of BF is partially supported by
RFBR initiative interdisciplinary project grant 09-02-12446-ofi-m,
by RFBR-CNRS grant 09-02-93106, RFBR grants 08-01-00720-a,
NSh-3472.2008.2 and 07-01-92214-CNRSL-a.
Research of EF was partially supported by the Russian President Grant MK-281.2009.1,
the RFBR Grants 09-01-00058, 07-02-00799 and NSh-3472.2008.2, by Pierre
Deligne fund based on his 2004 Balzan prize in mathematics and by
Alexander von Humboldt Fellowship.
Research of MJ is supported by the Grant-in-Aid for Scientific
Research B-20340027.
Research of EM is partially
supported by NSF grant DMS-0900984. Most of the present work has been carried out during the visits
of BF, EF and EM to Kyoto University. They wish to thank the University for hospitality.

\end{document}